\documentclass[12pt]{amsart}

\usepackage{amsmath,amsfonts,amsthm,amssymb,color}
\usepackage[latin1]{inputenc}
\usepackage[all,cmtip]{xy}

\numberwithin{equation}{section}

\newtheorem{ftheorem}{Th\'eor\`eme}
\newtheorem{fcor}{Corollaire}
\newtheorem{fconj}{Conjecture}
\newtheorem{fdef}{D\'efinition}
\newtheorem{prop}{Proposition}

\newcommand{\pa}{\partial}
\newcommand{\eps}{\varepsilon}

\newcommand{\n}{\nabla}

\newcommand{\N}{\mathbb{N}}
\newcommand{\M}{\mathfrak{M}}
\newcommand{\R}{\mathbb{R}}

\newcommand{\cC}{\mathcal{C}}

\renewcommand{\epsilon}{\varepsilon}

\begin{document}

\title [La g\'eom\'etrie de Bakry-\'Emery et l'\'ecart fondamental]{La g\'eom\'etrie 
de Bakry-\'Emery et l'\'ecart fondamental}
\author[J. Rowlett]{Julie Rowlett}\address{ Hausdorff Center for Mathematics\\ Universit\"at Bonn}
\date{}




\begin{abstract}
Cet article est une pr\'esentation rapide, d'une part de r\'esultats de l'auteur et 
Z. Lu \cite{lr}, et
d'autre part, de la r\'esolution de la conjecture de l'\'ecart fondamental 
 par Andrews et Clutterbuck \cite{ac}.
Nous commen\c cons  par rappeler ce qu'est la g\'eom\'etrie de Bakry-\'Emery, nous poursuivons 
en montrant les
liens entre valeurs propres du laplacien de Dirichlet et de Neumann.  Nous d\'emontrons 
ensuite un
rapport entre l'\'ecart fondamental et la g\'eom\'etrie de Bakry-\'Emery, puis nous pr\'esentons
les id\'ees principales de la preuve de la conjecture de l'\'ecart fondamental de \cite{ac}.  
Nous concluons par des r\'esultats pour l'\'ecart des triangles et des simplexes.
\end{abstract}

\maketitle

\section{La g\'eom\'etrie de Bakry-\'Emery}
Soit $(M, g)$ une vari\'et\'e riemannienne (avec ou sans bord); Bakry et \'Emery ont introduit
une g\'eom\'etrie qu'ils ont utilis\'ee pour \'etudier les processus de diffusion \cite{be}.
Une \em vari\'et\'e Bakry-\'Emery \em est un triplet $(M, g, \phi)$, o\`u la fonction 
$\phi \in \cC^{\infty} (M)$. La mesure sur $M$ est la mesure \`a poids $e^{-\phi} dV_g$, o\`u $dV_g$ 
est la mesure associ\'ee \`a la m\'etrique $g$. Le laplacien de Bakry-\'Emery est donn\'e par
$$\Delta_{\phi} = \Delta_g - \nabla \phi \cdot \nabla.$$
La courbure de Bakry-\'Emery-Ricci est\footnote{Dans Lott \cite{lott}, elle est appel\'ee
$\infty$-courbure-Bakry-\'Emery-Ricci.}
$${\rm Ric}_{\infty} = {\rm Ric} + {\rm Hess} (\phi).$$
On s'int\'eresse \`a la g\'eom\'etrie de Bakry-\'Emery pour g\'en\'eraliser 
la g\'eom\'etrie diff\'erentielle
aux vari\'et\'es singuli\`eres; voir Sturm \cite{st}, Wei-Wylie \cite{ww}, et 
Lott \cite{lott}.  Il est possible
de g\'en\'eraliser la notion de courbure de Ricci aux vari\'et\'es singuli\`eres 
qui sont la limite de Gromov-Hausdorff point\'ee de vari\'et\'es riemanniennes lisses \`a courbure de 
Ricci born\'ee inf\'erieurement.
Ces limites sont des espaces m\'etriques-mesur\'es, et sont aussi \'etudi\'ees 
dans le transport optimal; voir Villani \cite{vill}.

Le r\'esultat suivant montre que pour une vari\'et\'e de Bakry-\'Emery de dimension $n$ donn\'ee, 
il existe une famille \`a un param\`etre (positif) de domaines de dimension $n+1$, 
qui s'effondrent sur la vari\'et\'e quand le param\`etre tend vers $0$ de sorte que 
la famille \`a un param\`etre de  valeurs propres associ\'ees convergent vers celles de 
la vari\'et\'e. 

\begin{ftheorem}[Lu-R.]
Soit $(M, g, \phi)$ une vari\'et\'e de Bakry-\'Emery (avec ou sans bord).  Soit
$$M_{\eps} := \{ (x, y) \in M \times \R^+ |  \quad 0 \leq y \leq \eps e^{-\phi(x)} \} 
\subset M \times \R^+.$$
Soient $\{\mu_k \}_{k=0} ^{\infty}$ les valeurs propres du laplacien de Bakry-\'Emery de $M$; 
lorsque
$\pa M \neq \emptyset$ on consid\`ere la condition Neumann au bord. Soient $\mu_{k, \eps}$ 
les valeurs
propres (de Neumann si $\pa M \neq \emptyset$) du laplacien
$$\tilde{\Delta} := \Delta_g + \pa_y ^2,$$
sur $M_\eps$, o\`u $\Delta_g$ est le laplacien pour la m\'etrique $g$ sur $M$.  
Alors, $$\lim_{\eps \to 0} \mu_{k, \eps} = \mu_k \quad \forall \quad k \in \N.$$
\end{ftheorem}

Une cons\'equence imm\'ediate de ce th\'eor\`eme est le 

\begin{fcor}[Lu-R.]
Soit $\Omega$ un domaine de $\R^n$, et soit $\phi_1$ la premi\`ere fonction propre du laplacien 
euclidien sur $\Omega$.  Soit $$\Omega_\eps :=\{(x,y)\in\mathbb R^{n+1}\mid x\in \Omega,0\leq y\leq
\eps\phi_1(x)^2\}.$$
Soient $\{ \lambda_k\}_{k=1} ^{\infty}$ les valeurs propres de Dirichlet de $\Omega$, et soient
$\{ \mu_{k,\eps} \}_{k=0} ^{\infty}$ les valeurs propres de Neumann de $\Omega_\eps$.  Alors,
$$\lim_{\eps \to 0} \mu_{k-1,\eps} = \lambda_k-\lambda_1, \quad
\forall \quad k \in \N, \, k \geq 1.$$
\end{fcor}

Rappelons que pour un domaine $\Omega $ de $\R^n$, \em l'\'ecart fondamental de $\Omega$ \em 
est la  diff\'erence entre les deux premi\`eres valeurs propres de Dirichlet, $\lambda_2 - \lambda_1$.   
Le corollaire a des implications int\'eressantes pour l'\'ecart fondamental; pour $k=2$ on a 
$$\lim_{\eps \to 0} \mu_{1, \eps} = \lambda_2 - \lambda_1.$$

\subsection{Techniques utiles}

Rappelons les principes variationnels classiques.  Notre convention pour le laplacien 
associ\'e \`a une m\'etrique riemannienne $g$ est
$$\Delta = \frac{1}{\sqrt{\det(g)}} \sum_{i,j} \pa_i g^{ij} \sqrt{\det(g)} \pa_j.$$
Le laplacien euclidien est donc
$$\Delta =  \sum_{j=1} ^n \frac{\pa^2 }{\pa x_j ^2}.$$
Les valeurs propres de Dirichlet (de Neumann) sont les nombres r\'eels $\lambda$ pour 
lesquels il existe une fonction, alors dite propre, $u \in \cC^{\infty} (\Omega)$ telle que
$$ -\Delta u = \lambda u \textrm{ et } \left. u \right|_{\pa \Omega} = 0,
\textrm{ (Neumann : } \left. \frac{\pa u}{\pa n} \right|_{\pa \Omega} =0),$$
o\`u $n$ est le champ de vecteur normal de $\pa \Omega$.  Les valeurs propres de 
Dirichlet sont not\'ees $\lambda$ et indic\'ees par $\N_{\geq 1 }$. Les valeurs propres 
de Neumann sont not\'ees 
$\mu$ et indic\'ees par $\N$. Les valeurs propres de Dirichlet et de Neumann satisfont 
les principes variationnels (voir par exemple Chavel \cite{chavel} et
Courant-Hilbert \cite{cour-hill})
$$ \lambda_1 = \inf_{f \in \cC^1 (\Omega)} \left\{ \left. 
\frac{ \int_{\Omega} |\nabla f |^2 }{\int_{\Omega} f^2 }  \,\right|  \, \,  \left. 
f\right|_{\pa \Omega} = 0, \, f \not\equiv 0 \right\},$$
$$\mu_0 = \inf_{f \in \cC^1 (\Omega)} \left\{ \left. 
\frac{ \int_{\Omega} |\nabla f |^2 }{\int_{\Omega} f^2 } \,\right| \, \, 
f \not\equiv 0 \right\},$$
et pour $k > 1$, $j > 0$,
$$\lambda_k = \inf_{f \in \cC^1 (\Omega)} \left\{ \left. 
\frac{ \int_{\Omega} |\nabla f |^2 }{\int_{\Omega} f^2 } \, \right| \, \left. \,  
f\right|_{\pa \Omega} = 0, \, f \not\equiv 0 = \int_{\Omega} f \phi_j, \, 0<  j < k 
\right\},$$
$$\mu_j = \inf_{f \in \cC^1 (\Omega)} \left\{   \left. 
\frac{ \int_{\Omega} |\nabla f |^2 }{\int_{\Omega} f^2 } \, \right| \, \,  
f \not\equiv 0 = \int_{\Omega} f \varphi_l, \, 0 \leq l < j \right\},$$
o\`u $\phi_j$ et $\varphi_l$ sont les fonctions propres associ\'ees \`a $\lambda_j$ et 
\`a $\mu_j$, respectivement. Si on n'impose aucune condition au bord, la condition de 
Neumann est alors imm\'ediatement satisfaite par la fonction r\'ealisant l'infimum.

Les valeurs propres satisfont aussi les principes ``min-max'' suivants:
$$\lambda_k = \inf \left\{ \sup \left . \left\{ \left . \frac{\int_\Omega|\nabla f|^2}{\int_\Omega f^2}
\right|  \, f \in L \right\} \right|  \, \textrm{dim}(L) = k, \, f|_{\pa \Omega} = 0 \, \forall 
\, f \in L \right\},$$
$$\mu_k = \inf \left\{ \sup \left . \left\{ \left . \frac{\int_\Omega|\nabla f|^2}{\int_\Omega f^2}
\right|  \, f \in L \right\} \right|  \,  \textrm{dim}(L) = k \right\},$$
o\`u $L \subset H^1(\Omega)$.  Les principes variationnels et min-max sont identiques 
pour un op\'erateur Schr\"odinger $\Delta + V$, o\`u le potentiel $V$ est une 
fonction lisse.

Rappelons la proposition suivante \cite{lr}: 
\begin{prop}\label{prop1}  Soit $\Omega \subset \R^n$ un domaine de $\R^n$ \`a bord lisse et 
soit une base orthonormale de fonctions propres $\{ \phi_k \}_{k=1} ^{\infty}$ respectivement 
associ\'ees aux valeurs propres Dirichlet $\{ \lambda_k \}_{k=1} ^{\infty}$.
Alors, pour tout $k \in \N$, $\psi_k = \frac{\phi_k}{\phi_1}$ est lisse jusqu'au bord et satisfait:
\begin{equation}\label{eq:1}
\Delta \psi_k + 2\nabla\log \phi_1\nabla\psi_k= - (\lambda_k - \lambda_1) \psi_k.
\end{equation}
\begin{equation}\label{eq:2} \left.\frac{\pa \psi_k}{\pa n}\right|_{\pa\Omega}=0.
\end{equation}
De plus, $\n \log \phi_1 \n \psi_k$ est lisse jusqu'au bord.  Si le bord est seulement lisse par 
morceaux, l'\'egalit\'e (\ref{eq:1}) est encore satisfaite, $\n \log \phi_1 \n \psi_k$ est 
lisse jusqu'aux parties lisses de $\pa \Omega$, et de plus l'\'egalit\'e (\ref{eq:2}) est 
v\'erifi\'ee  sur chaque composante lisse de $\pa \Omega$.
\end{prop}

Une cons\'equence imm\'ediate de la proposition pr\'ec\'edente
(pour $k=2$ on retouve un r\'esultat de Ma-Liu \cite{maliu})
est

\begin{prop}[Lu-R./Ma-Liu] \label{prop2}
Soit $\Omega \subset \R^n$ un domaine born\'e de $\R^n$ et soit  $\{ \phi_k \}_{k=1} ^{\infty}$ une
base orthonormale de fonctions propres respectivement associ\'ees aux valeurs propres
de Dirichlet $\{ \lambda_k \}_{k=1} ^{\infty}$ du laplacien eucildien. Soient $\{\mu_k\}_{k=0} ^{\infty}$ 
les valeurs propres de Neumann du laplacien de Bakry-\'Emery pour la fonction poid $-2 \log \phi_1$.  
Alors, $$\lambda_k - \lambda_1 = \mu_{k-1} \qquad \forall \quad k \in \N, \, k \geq 1.$$
\end{prop}

Motiv\'e par cette proposition, nous avons d\'emontr\'es les principes variationnels 
pour le laplacien de Bakry-\'Emery;
pour $k=2$ et $M \subset \R^n$, c'est le corollaire 1.3 de Kirsch-Simon \cite{kir-sim}.

\begin{prop}[Lu-R./Kirsch-Simon]\label{prop3}  Soit $(M, g, \phi)$ une vari\'et\'e de Bakry-\'Emery 
(avec ou sans bord).  Les valeurs propres du laplacien de Bakry-\'Emery satisfont:
\\
pour la condition Dirichlet si $\pa M \neq \emptyset$:
$$
\lambda_1 = \inf_{\varphi \in \cC^1(M)} \left \{ \left . 
\frac{\int_M|\nabla\varphi|^2e^{- \phi}}{\int_M\varphi^2 e^{-\phi} }
\, \right| \, \, \varphi \not\equiv 0, \, \varphi |_{\pa M} = 0 \right\};
$$
pour la condition Neumann si $\pa M = \emptyset$:
$$
\mu_0 = \inf_{\varphi \in \cC^1(M)} \left \{ \left . 
\frac{\int_M|\nabla\varphi|^2e^{- \phi}}{\int_M\varphi^2 e^{-\phi} }
\, \right| \, \, \varphi \not\equiv 0, \right\}.
$$
Pour $k \geq 2$,
$$\lambda_k =  \inf_{\varphi \in \cC^1(M)}  \left \{ \left . 
\frac{\int_M|\nabla\varphi|^2 e^{-\phi} }{\int_M\varphi^2 e^{-\phi}}
\, \right|  \, \,  \varphi \not\equiv 0 = \int_{M} \varphi \varphi_j e^{-\phi}, \, \, 1 \leq j < k, \,
\varphi |_{\pa M} = 0 \right\},
$$
$$\mu_k =  \inf_{\varphi \in \cC^1(M)}  \left \{ \left . 
\frac{\int_M|\nabla\varphi|^2 e^{-\phi} }{\int_M\varphi^2 e^{-\phi}}
\, \right|  \, \,  \varphi \not\equiv 0 = \int_{M} \varphi \varphi_j e^{-\phi}, \, 0 \leq j < k \right\},
$$
o\`u $\varphi_j$ est la fonction minimisante quand $k=j$.
\end{prop}

Si $M \subset \R^n$ avec valeurs propres de Dirichlet $\{ \lambda_k \}_{k=1} ^{\infty}$ et fonctions
propres orthonormales associ\'ees $\{ \phi_k \}_{k=1} ^{\infty}$, et si la fonction poid est 
$\phi = - 2 \log \phi_1$, le principe variationnel pour $(M, g_{eucl} , \phi)$ est alors
$$
\lambda_k - \lambda_1 = \inf_{\varphi \in \cC^1(\Omega)} \left \{ \left .
\frac{\int_\Omega|\nabla\varphi|^2\phi_1^2}{\int_\Omega\varphi^2\phi_1^2} \, \right| \, \, 
\varphi \not\equiv 0 \right\},$$
et pour $k \geq 2$,
$$\lambda_k - \lambda_1 =  \inf_{\varphi \in \cC^1(\Omega)}  \left \{ \left .
\frac{\int_\Omega|\nabla\varphi|^2\phi_1^2}{\int_\Omega\varphi^2\phi_1^2} \, \right|  \, \,
\varphi \not\equiv 0 = \int_{\Omega} \varphi \varphi_j \phi_1^2, \, 1 \leq j < k \right\}.$$

Naturellement, si il y a un principe variationnel, il y a aussi un principe min-max.

\begin{prop}[Lu-R.]\label{prop3}  Soit $(M, g, \phi)$ une vari\'et\'e de Bakry-\'Emery 
(avec ou sans bord). Alors les valeurs propres du laplacien de Bakry-\'Emery satisfont:
pour la condition Dirichlet si $\pa M \neq \emptyset$
$$
\lambda_k = \inf \left\{ \sup \left . \left\{ \left . 
\frac{\int_M|\nabla\varphi|^2 e^{-\phi}}{\int_M\varphi^2 e^{-\phi}}
\right|  \, \varphi \in L \right\} \right|  \, \textrm{dim}(L) = k, \, f|_{\pa \Omega} = 0 \, \forall
\, f \in L \right\},
$$
et pour la condition Neumann si $\pa M = \emptyset$:
$$
\mu_k = \inf \left\{ \sup \left . \left\{ \left . 
\frac{\int_M |\nabla \varphi |^2 e^{-\phi}}{\int_M \varphi^2 e^{-\phi}}
\right|  \, \varphi \in L \right\} \right|  \, \textrm{dim}(L) = k \right\}, 
$$
o\`u $L \subset H^1(M, e^{-\phi} dV_g)$.\\
\end{prop}

Une estimation utile est donn\'ee par la proposition suivante:

\begin{prop}[Lu-R.]\label{prop4} Soient $k \geq 1$, et $\xi_1,\cdots,\xi_{k}$ orthogonales 
(et non triviales) dans $(M, g, \phi)$, c'est-\`a-dire que $\xi_i\not\equiv 0$, $1 \leq i \leq k$ et 
\[
\int_M \xi_i\xi_j e^{-\phi}=0, \quad i \neq j.
\]
Alors, si $\pa M = \emptyset$, ou si $\pa M \neq \emptyset$ et on consid\`ere alors la condition de 
Neumann au bord, les valeurs propres du laplacien de Bakry-\'Emery $\{ \mu_j \}$ satisfont   
\[
\sum_{j=0}^{k} \mu_j\leq\sum_{j=1}^{k}\frac{\int_M |\n\xi_j|^2 e^{-\phi}}{\int_M |\xi_j|^2 e^{-\phi}}.
\]
\end{prop}

\section{L'\'ecart fondamental}
La fonction \em \'ecart $\xi$ \em sur l'espace des domaines de $\R^n$ est d\'efinie comme suit:
$$\xi: M \to \R, \quad \xi(M) = d^2 (\lambda_2 - \lambda_1),$$
o\`u $d$ est le diam\`etre de $M$ et $\lambda_1 < \lambda_2$ sont les premi\`eres valeurs propres du 
laplacien euclidien avec la condition de Dirichlet au bord.  
Le \em probl\`eme de l'\'ecart \em consiste \`a estimer la fonction \'ecart sur l'espace des
domaines \em convexes. \em  Lorsque les domaines ne sont pas convexes, en consid\'erant des domaines qui
ont la forme de deux sph\`eres separ\'ees par une longue tube, l'\'ecart tend vers zero.  Si le domaine
est normalis\'e de sorte \`a avoir diam\`etre un, la fonction \'ecart est alors \'egale \`a 
\em l'\'ecart fondamental :
\em
$$\lambda_2 - \lambda_1.$$
Van den Berg \cite{mvdb} a conjectur\'e que la fonction \'ecart (sur les domaines convexes) est 
born\'ee inf\'erieurement par une constante; il est donc naturel de conjecturer que la constante 
est $3 \pi^2$, l'\'ecart de l'intervalle $[0,1]$.  Par exemple, on peut calculer les valeurs propres 
de Dirichlet d'un rectangle $R \cong [0,a] \times [0,b]$, o\`u $a \geq b$.  Les valeurs 
propres sont alors $$\lambda_{j,k} = \frac{\pi^2 j^2}{a^2} + \frac{\pi^2 k^2}{b^2},$$
et l'\'ecart fondamental est 
$$\frac{3 \pi^2}{a^2}.$$
La fonction \'ecart est 
$$\frac{3 \pi^2 (a^2 + b^2)}{a^2}.$$
On voit dans cet exemple que la fonction \'ecart sur les rectangles est maximis\'ee par 
le carr\'e et tend vers son infimum $3 \pi^2$ quand un rectangle s'\'ecrase sur un intervalle.  

Les premiers r\'esultats estimant l'\'ecart utilisent des
estimations de gradient dans l'esprit de Li-Yau \cite{li-yau};
Singer-Wong-Yau-Yau \cite{swyy} ont montr\'e que l'\'ecart satisfait
$$\xi \geq \frac{\pi^2}{4}.$$
En raffinant p\'eniblement les m\^emes techniques, Yu-Zhong \cite{yz}  prouve l'estimation
$$\xi \geq \pi^2;$$
voir aussi Li-Treibergs \cite{lt}.  

Le corollaire 1 montre que si le diam\`etre de $M$ est un alors
$$\xi = \underset{\eps\to 0}{\lim}\,\mu_{1,\eps}.$$
En utilisant ce corollaire directement avec les estimations de gradient de  \cite{li-yau}, 
on obtient les r\'esultats de  \cite{swyy} et \cite{yz} avec des 
preuves beaucoup plus courtes.
Nous avons remarqu\'e \cite{lr} que le Hessien du logarithme de la premi\`ere fonction propre 
joue le r\^ole de la courbure Ricci pour les estimations de gradient de Li-Yau. On peut obtenir 
tous les r\'esultats d'estimations de gradient pour la g\'eom\'etrie de Bakry-\'Emery; 
c'est un projet int\'eressant que d'appliquer de tels r\'esultats aux espaces m\'etriques-mesur\'es 
\cite{st}.  Dans tous les cas, l'observation que \em le Hessien du logarithme de la premi\`ere 
fonction propre joue le r\^ole de la courbure Ricci \em a \'et\'e une des cl\'es de la preuve 
de la conjecture de l'\'ecart fondamental.

L'autre cl\'e est l'\'ecart associ\'e \`a un op\'erateur Schr\"odinger \`a potentiel convexe en dimension
une; en dimension une, ce probl\`eme a \'et\'e resolu par Lavine \cite[Theorem 3.1]{lav}.

\begin{ftheorem}[Lavine] Soit $V$ une fonction convexe sur $[0, R]$ et soient $\lambda_1$, 
$\lambda_2$ les deux premi\`eres valeurs propres de Dirichlet (respectivement Neumann) 
de l'op\'erateur Schr\"odinger $-d^2/dx^2 + V$ sur $[0, R]$.  Alors
$$\lambda_2 - \lambda_1 \geq \Gamma_0,$$
o\`u $\Gamma_0$ est l'\'ecart fondamental pour $V =  constante$ pour l'op\'erateur de Dirichlet
(respectivement de Neumann), et il y a \'egalit\'e si et seulement si $V =  constante$.
Alors, pour l'op\'erateur de Dirichlet, $$\lambda_2 - \lambda_1 \geq \frac{3\pi^2}{R^2},$$
et pour l'op\'erateur de Neumann $$\lambda_2 - \lambda_1 \geq \frac{\pi^2}{R^2}.$$
\end{ftheorem}

Pour utiliser le th\'eor\`eme de Lavine, Andrews et Clutterbuck introduisent les \em module de 
continuit\'e \em et \em module de convexit\'e/concavit\'e.  \em
\begin{fdef}[Andrews-Clutterbuck]
Soient une fonction $\eta: \R^+ \to \R^+$ et une fonction $f: \Omega \to \R$, o\`u $\Omega \subset \R^n$
est un domaine.  La fonction $\eta$ est un \em module de continuit\'e \em pour la fonction $f$ si
$$| f(y) - f(x) | \leq 2 \eta \left( \frac{|y-x|}{2} \right).$$
\end{fdef}

Pour d\'efinir le \em module de convexit\'e et concavit\'e \em il faut d'abord d\'efinir le \em 
module d'expansion \em et \em module de contraction.  \em  
\begin{fdef}[Andrews-Clutterbuck]
Soit une fonction $\omega: \R^+ \to \R$ et soit $X$ un champs de vecteur d\'efini sur un domaine
$\Omega \subset \R^n$.  La fonction $\omega$ est un \em module d'expansion \em pour $X$ si
$$(X(y) - X(x)) \cdot \frac{y-x}{|y-x|} \geq 2 \omega \left( \frac{|y-x|}{2} \right),
\quad \forall \,  x, y \in \Omega,  \, y \neq x.$$
La fonction $\omega$ est un \em module de contraction \em pour $X$ si $-\omega$ est un module d'expansion pour $-X$.  
\end{fdef}

\begin{fdef}[Andrews-Clutterbuck]
Soit une fonction $\omega: \R^+ \to \R$ et soit $f$ une fonction semi-convexe sur un domaine $\Omega \subset \R^n$.  La fonction $\omega$ est un \em module de convexit\'e \em pour $f$ si $\omega$ est un module d'expansion pour le champs de vecteur donn\'e par le gradient de la fonction $f$, $\nabla f$.  La fonction $\omega$ est un \em module de concavit\'e \em pour $f$ si elle est un module de contraction pour $\nabla f$.  
\end{fdef}

Ensuite, Andrews et Clutterbuck ont utilis\'e ces id\'ees pour raffiner le r\'esultat de 
Brascamp-Lieb \cite{bl}:  Hess $\log \phi_1 \leq 0$.  Andrews et Clutterbuck d\'emontrent 
que $(\log \tilde{\phi}_1)'$ est un module de concavit\'e pour $\log \phi_1$, o\`u $\tilde{\phi}_1$ 
est la fonction propre pour un op\'erateur Schr\"odinger en dimension une.  Ce r\'esultat signifie 
que $\phi_1$ est ``davantage log-concave'' que la premi\`ere fonction propre du probl\`eme en 
dimension une.  Puisque le probl\`eme a \'et\'e resolu en dimension une, Andrews et Clutterbuck 
obtiennent l'estimation pour l'\'ecart en dimension $n$  \cite[Theorem 1.3]{ac} suivante:

\begin{ftheorem}[Andrews-Clutterbuck]
Soit $\Omega \subset \R^n$ un domaine convexe  de diam\`etre $R$.  Soient $\lambda_1$ et 
$\lambda_2$ les premi\`eres valeurs propres Dirichlet de $\Delta + V$ sur $\Omega$. Soit une 
fonction $\tilde{V} \in \cC^1 ( [0, R])$, et soient $E_1$ et $E_2$ les premi\`eres valeurs propres 
Dirichlet de $- d^2/dx^2 + \tilde{V}$ sur $[0, R]$.  Si $\tilde{V}'$ est un module de convexit\'e 
pour $V$, alors 
$$\lambda_2 - \lambda_1 \geq E_2 - E_1.$$
\end{ftheorem}

Enfin, ils ont montr\'e la conjecture de l'\'ecart fondamental.
\begin{fcor}[Andrews-Clutterbuck]
Si $V$ est convexe alors
$$\lambda_2 - \lambda_1 \geq \frac{3 \pi^2}{R^2}.$$
\end{fcor}

Si on consid\`ere l'\'ecart sur l'espace des modules du triangle, il se comporte 
totalement diff\'eremment.  

\section{L'\'ecart de simplexes}
Dans la derni\`ere section  de cet article, nous consid\'erons un probl\`eme plus concret.  
Rappelons qu'un $n$-simplexe $X$ est un ensemble de $n+1$ vecteurs $\{v_0,\cdots,v_n\}$ de 
$\mathbb R^n$ tels que $v_1-v_0,\cdots,v_n-v_0$ sont lin\'eairement ind\'ependants. Le domaine convexe
$$
\left\{\left . \sum_{j=0}^n t_j v_j \right| \sum_{j=0} ^n t_j=1, \,  t_j\geq 0
\textrm{ for } 0 \leq j \leq n \right\}
$$
est l'enveloppe convexe de $\{v_0,\cdots,v_n\}$.  Ce domaine est born\'e, \`a bord lisse par 
morceaux. Par simplicit\'e, nous \'ecrivons le simplexe $X$ pour signifier \`a la fois 
$\{v_0,\cdots,v_n\}$ et l'enveloppe convexe de $\{v_0,\cdots,v_n\}$.  Si $n=2$, un simplexe est 
simplement un triangle. L'espace des modules du $n$-simplexe est parametr\'e par l'ensemble 
 des $n$-simplexes de diam\`etre un. En contraste du th\'eor\`eme de Andrews et Clutterbuck, 
la fonction \'ecart sur l'espace des modules du $n$-simplexe \em n'est pas \em minimis\'e par l
es simplexes qui s'effondrent en dimension inf\'erieure.

\begin{ftheorem}[Lu-R.]\label{th:simplex}
Soit $Y$ un $(n-1)$-simplexe, o\`u $n \geq 2$.  Soit $\{X_j \}_{j \in
\N}$ une suite de $n$-simplexes, telle que chacun est un graphe sur $Y$ et la hauteur de
 $X_j$ sur $Y$ tend vers zero quand $j \to \infty$. Alors,
$$\xi(X_j) \to \infty, \quad j \to \infty.$$
Il existe alors une constante $C>0$ d\'ependant seulement de $n$ et de $Y$ telle que
$$\xi(X_j)\geq C h(X_j)^{-4/3},$$ o\`u $h(X_j)$ est la hauteur de $X_j$ sur $Y$.
\end{ftheorem}

Dans le cas des triangles ($n=2$), le th\'eor\`eme implique alors le corollaire 
suivant:\footnote{Une preuve
ind\'ependante de ce th\'eor\`eme se trouve dans \cite{lr1}.}

\begin{fcor}[Friedlander-Solomyak] \label{compactness}
Soit $P$ l'espace des modules du triangle.  Alors, la fonction \'ecart $\xi: P \to \R$ est propre.
\end{fcor}

Le corollaire implique l'existence d'un triangle minimisant l'\'ecart; en collaboration avec 
T. Betcke \cite{lrb}, on montre la conjecture d'Antunes-Freitas \cite{freitas}:

\begin{ftheorem}[Betcke-Lu-Rowlett] \label{th:eqtri} Soit $T$ un triangle.  Alors,
$$\xi(T) \geq \frac{64 \pi^2}{9},$$
et il y a \'egalit\'e si et seulement si $T$ est \'equilat\'eral.
\end{ftheorem}

Il est alors naturel de proposer la g\'en\'eralisation suivante:
\begin{fconj}  Soit $\M_n$ l'espace des modules du $n$-simplexe, $n \geq 2$.  Alors la fonction \'ecart
$\xi: M_n \to \R$ est propre, et le simplexe d\'efini par $p_0, p_1, \ldots, p_n \in \R^n$ de sorte que
$$|p_i - p_j| = 1 \textrm{ for } 0 \leq i \neq j \leq n$$
minimise uniquement la fonction \'ecart sur $\M_n$.\\
\end{fconj}


\begin{thebibliography}{99}

\bibitem[A-C10]{ac} B. Andrews et J. Clutterbuck, \em Proof of the fundamental gap conjecture,
\em arXiv 1006.1686, (2010).

\bibitem[A-F08]{freitas} P. Antunes et P. Freitas, \em A numerical study of the spectral gap,
\em J. Phys. A 41, no. 5, (2008), 055201, 19 p.
	
\bibitem[B-E83]{be} D. Bakry et M. \'Emery, \em Diffusions hypercontractives, \em
S\'eminaire de probabilites, XIX, 1983/84, Lecture Notes in Math., 1123, Springer, Berlin, (1985) 177--206.

\bibitem[B-L-R]{lrb} T. Betcke, Z. Lu, et J. Rowlett, \em The fundamental gap of triangles, 
\em en pr\'eparation.

\bibitem[B-L76]{bl} H. J. Brascamp et E. H. Lieb, \em On extensions of the Brunn-Minkowski and 
Pr\'ekopa-Leindler
   theorems, including inequalities for log concave functions, and with an
   application to the diffusion equation, \em J. Functional Analysis, 22, no. 4, (1976), 366--389.

\bibitem[C74]{chavel} I. Chavel, \em Eigenvalues in Riemannian geometry, \em
Pure and Applied Mathematics, vol. 115, Academic Press Inc., Orlando, FL, (1984).

\bibitem[C-H53]{cour-hill} R. Courant et D. Hilbert, \em Methods of mathematical physics. Vol. I,
\em Interscience Publishers, Inc., New York, N.Y., (1953).

\bibitem[F-S07]{strip} L. Friedlander et M. Solomyak, \em On the Spectrum of the Dirichlet Laplacian in a
Narrow Strip, \em arXiv math.SP 0705.4058v1, (2007).

\bibitem[K-S87]{kir-sim} W. Kirsch et B. Simon, \em Comparison theorems for the gap of 
Schr\"odinger operators,
\em J. Funct. Anal., 75, no. 2, (1987) 396--410.

\bibitem[L94]{lav} R. Lavine, \em The eigenvalue gap for one-dimensional convex potentials,
\em Proc. Amer. Math. Soc., 121, no. 3, (1984), 815--821.

\bibitem[L-T91]{lt} P. Li et A. Treibergs, \em  Applications of eigenvalue techniques to geometry, 
\em Contemporary geometry, Univ. Ser. Math., Plenum, New York, (1991), 21--52. 

\bibitem[L-Y79]{li-yau} P. Li et S. T. Yau, \em Estimates of eigenvalues of a compact Riemannian manifold,
\em Geometry of the Laplace operator,  Proc. Sympos. Pure Math. XXXVI, Amer. Math. Soc.,
Providence, R. I., (1980) 205--239.

	
\bibitem[L03]{lott}  J. Lott, \em Some geometric properties of the Bakry-\'Emery-Ricci tensor,
\em
Comment. Math. Helv., 78, no. 4, (2003), 865--883.

\bibitem[L-R09]{lr} Z. Lu et J. Rowlett, \em The fundamental gap, \em preprint (2009), arXiv 1003.0191v1.

\bibitem[L-R08]{lr1} Z. Lu et J. Rowlett, \em The Fundamental Gap Conjecture on Polygonal Domains,
\em arXiv:0810.4937, (2008).

\bibitem[M-L09]{maliu} L. Ma et B. Liu, \em Convex eigenfunction of a drifting Laplacian operator and
the fundamental gap, \em  Pacific J. Math. 240, no. 2, (2009), 343--361.
	
\bibitem[SWYY85]{swyy} I. M. Singer, B. Wong, S. T. Yau, et S. S. T. Yau, \em An estimate of the gap of
the first two eigenvalues in the Schr\"odinger operator, \em Ann. Scuola Norm. Sup. Pisa Cl. Sci. (4), 12,
no. 2, (1985), 319--333.

\bibitem[S06]{st} K. T. Sturm, \em On the geometry of metric measure spaces I, II, \em Acta Math. 196,
no. 1, (2006), 65--131 and 133--177.

\bibitem[V83]{mvdb} M. van den Berg, \em On the condensation of the free Boson gas and the spectrum of
the Laplacian, \em J. Stat. Phys. 31, (1983), 623--637.

\bibitem[V04]{vill} C. Villani, \em Topics in optimal transportation, \em
Graduate Studies in Mathematics, 58. American Mathematical Society, Providence, RI, (2003).

\bibitem[WW07]{ww} G. Wei et W. Wylie, \em Comparison geometry for the Bakry-Emery Ricci tensor, \em
J. Differential Geom. 83, no. 2, (2009), 377--405.

\bibitem[Y-Z86]{yz} Q. H. Yu et J. Q. Zhong, \em Lower bounds of the gap between the first and second
eigenvalues of the Schr\"odinger operator, \em Trans. Amer. Math. Soc., 294, no. 1, (1986), 341--349.

\end{thebibliography}
\end{document}